\documentclass[reqno]{amsart}
\usepackage{amsmath}%
\usepackage{amsfonts}%
\usepackage{amssymb}%
\usepackage{graphicx}
\usepackage{setspace}
\usepackage{microtype}

\theoremstyle{plain}

\numberwithin{equation}{section}
\begin{document}
\onehalfspacing
\raggedbottom
\title{Constructing Cartesian Splines}
\author{H.R.N.~van~Erp}
%\address[Structural Hydraulic Engineering and Probabilistic Design, TU DElft, The Netherlands]
%{Author OneTwo common address, line 1 \newline%
%\indent Author OneTwo common address, line 2}%
%\email[A. One]{author-one@autherone-inst.de}%
%\urladdr{http://www.authorone.uni-aone.de}
\author{R.O.~Linger}
%\curraddr[A.~Two]{Author Two current address, line 1\newline%
%\indent Author Two current address, line 2}%
%\email[A.~Two]{author-two@authortwo-inst.hu}%
%\urladdr{http://www.authortwo.uni-atwo.hu}
\author{P.H.A.J.M.~van~Gelder}

%\author{&N.~van~Erp\\Structural Hydraulic and Probabilistic Design, TU Delft, Delft, The Netherlands\\
%&\and R.~Linger\\Clinical Epidemiology, UMCG, Groningen, The Netherlands\\
%&\and P.~van~Gelder\\Structural Hydraulic and Probabilistic Design, TU Delft, Delft, The Netherlands}

\begin{abstract}
\noindent We introduce here Cartesian splines or, for short, C-splines. C-splines are piecewise polynomials which are defined on adjacent Cartesian coordinate systems and are $C^{r}$ continuous throughout. The $C^{r}$ continuity is enforced by constraining the coefficients of the polynomial to lie in the null-space of some smoothness matrix $H$. The matrix-product of the null-space of the smoothness matrix $H$ and the original polynomial base results in a new base, the so-called C-spline base, which automatically enforces $C^{r}$ continuity throughout. In this article we give a derivation of this C-spline base as well as an algorithm to construct C-spline models.  
\end{abstract}

\maketitle

\section{Introduction}
We introduce here Cartesian splines or, for short, C-splines. C-splines are piecewise polynomials which are defined on adjacent Cartesian coordinate systems and are $C^{r}$ continuous throughout. The $C^{r}$ continuity is enforced by constraining the coefficients of the polynomial to lie in the null-space of some smoothness matrix $H$. The matrix-product of the null-space of the smoothness matrix $H$ and the original polynomial base results in a new base, the so-called C-spline base, which automatically enforces $C^{r}$ continuity throughout. The idea of using the null-space of some smoothness matrix $H$ has been taken from the B-spline literature, where piecewise polynomials are defined on  adjacent triangular Barycentric coordinate systems,~\cite{Awanou}. It turns out that C-spline bases have a particular simple form. This makes it possible to give an explicit formulation of general C-spline bases. In this article we will give a general outline how to enforce continuity constraints by way of the smoothness matrix $H$. We then show how these constraints lead us to the C-spline base. Then we will give the explicit algorithm for constructing a bivariate C-spline base and show how to use this base to construct a C-spline model.

\section{Piecewise Polynomials \label{S:2}}
We start with the bivariate Cartesian $x$,$y$-coordinate system. We partition this initial coordinate system with origin $O = \left(0,0\right)$ in two adjacent coordinate systems, each with its own origin, $O = \left(0,0\right)$ and $\tilde{O} = \left(0,0\right)$. The geometry in terms of $x$ and $y$ may be depicted as:  

\begin{figure}[!h]
	\centering
		\includegraphics{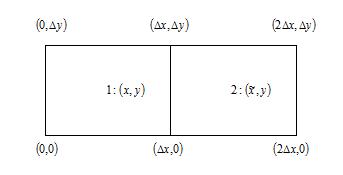}
	\caption{Geometry of the domain of two piecewise polynomials in terms of $x$ and $y$}
	\label{fig:stelsel1}
\end{figure}
\noindent where $\Delta x$ and $\Delta y$ are some constants. Likewise, the geometry in terms of $\tilde{x}$ and $y$ may be depicted as:

\begin{figure}[!h]
	\centering
		\includegraphics{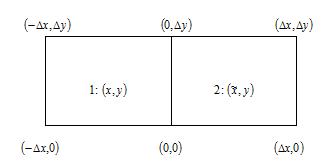}
	\caption{Geometry of the domain of two piecewise polynomials in terms of $\tilde{x}$ and $y$}
	\label{fig:stelsel2}
\end{figure}
\noindent where $\Delta x$ and $\Delta y$ are the same constants as used in Figure 1.

Now, we may define on both coordinate systems a polynomial of order $d$:
\begin{equation}
\label{E:0}
P_{d}\left(x,y\right) = \sum_{0\leq p + q\leq d}c_{pq} x^p y^q
\end{equation}
We start with the most simple case, that is, we set $d = 1$. The polynomial equations for both coordinate systems then become:
\begin{align}
\label{E:1}
z_{1}\left(x,y\right) = c_{11} + c_{12} x + c_{13} y, \quad\quad 0\leq x\leq \Delta x, \quad 0\leq y\leq \Delta y \notag \\
z_{2}\left(\tilde{x},y\right) = c_{21} + c_{22} \tilde{x} + c_{23} y, \quad\quad 0\leq \tilde{x}\leq \Delta x, \quad 0\leq y\leq \Delta y 
\end{align}
If we look at Figure 2, we see that 
\begin{equation}
\label{E:2}
\tilde{x} = x - \Delta x
\end{equation}
Combining~\eqref{E:1} and~\eqref{E:2} we get:
\begin{align}
\label{E:3}
z_{1}\left(x,y\right) &= c_{11} + c_{12} x + c_{13} y, \quad  &\;\;\;0\leq x\leq \Delta x,\quad \;\;0\leq y\leq \Delta y \notag \\
z_{2}\left(x,y\right) &= c_{21} + c_{22} \left(x - \Delta x\right) + c_{23} y, \quad &\Delta x\leq x\leq 2\Delta x,\quad 0\leq y\leq \Delta y 
\end{align}
Let 
\[
\mathbf{z} = \begin{pmatrix} z_{1}\left(x,y\right)&z_{2}\left(x,y\right)\end{pmatrix}^{T}
\]
be the outcome vector. Then~\eqref{E:3} may be rewritten as the matrix-vector product of the polynomial base
\begin{equation}
\label{E:5}
B = \begin{pmatrix} 
	 1&x&y&0&0&0\\
	 0&0&0&1&x - \Delta x&y
\end{pmatrix}
\end{equation}
and the coefficient vector
\begin{equation}
\label{E:6}
\mathbf{c} =  \begin{pmatrix} c_{11}&c_{12}&c_{13}&c_{21}&c_{22}&c_{23}\end{pmatrix}^{T}
\end{equation}
that is,
\[
\mathbf{z} = B \mathbf{c}
\]
Note that the $\left(x,y\right)$-values that fall in the first quadrant of Figure 1 are assigned to the first row of the polynomial base $B$, while $\left(x,y\right)$-values in the second quadrant are assigned to the second row. 

\section{Enforcing Zeroth Order Continuity \label{S:3}}
In order for the two polynomials~\eqref{E:3} to connect at the boundary, that is, in order to have $C^{0}$ continuity, we must have that
\begin{equation}
\label{E:8}
z_{1}\left(\Delta x,y\right) = z_{2}\left(\Delta x,y\right)
\end{equation}
for any $y$. Substituting~\eqref{E:3} in~\eqref{E:8}, we find  
\[
c_{11}+c_{12} \Delta x + c_{13} y = c_{21} + c_{23} y
\]
or, equivalently,
\begin{equation}
\label{E:8c}
c_{11}+c_{12} \Delta x + c_{13} y - c_{21} - c_{23} y = 0
\end{equation}
We have that~\eqref{E:8c} is a constraint on the $\mathbf{c}$ coefficients. The coefficients $\mathbf{c}$ must all lie in the null-space of the smoothness ``matrix'' $H$, where
\begin{equation}
\label{E:9}
H =  \begin{pmatrix} 1 &\Delta x &y &-1 &0 &-1\end{pmatrix}
\end{equation}
The null-space of $H$ is
\begin{equation}
\label{E:10}
H_{0} =  \begin{pmatrix} 
y &0 &1 &-y &-\Delta x\\
0 &0 &0 &0 &1 \\
0 &0 &0 &1 &0 \\
0 &0 &1 &0 &0 \\
0 &1 &0 &0 &0 \\
1 &0 &0 &0 &0 
\end{pmatrix}
\end{equation}
and it may be checked that
\[
H H_{0} = \mathbf{0} 
\]
where $\mathbf{0}$ is the $1\times5$ zero vector. It follows that the matrix product of $H$ with any linear combination of the columns in $H_{0}$ must give a zero value, that is,
\[
H H_{0} \mathbf{c}_{0} = 0
\]
where $\mathbf{c}_{0}$ is an arbitrary $5\times1$ vector. Stated differently, any linear combination of the columns of $H_{0}$ gives us an $6\times1$ vector that satisfies the constraint~\eqref{E:8} or, equivalently, constraint~\eqref{E:8c}. 

Now, if we take the matrix product of our original polynomial base, $B$, and the null-space of our smoothness matrix, $H_{0}$, we get the null-base $B_{0}$:
\begin{align}
\label{E:13}
B_{0} &= B H_{0}\notag\\
&= \begin{pmatrix} 
y &0 						&1 &0	&x -\Delta x \\
y &x -\Delta x	&1 &0 &0 
\end{pmatrix}
\end{align}
If we drop the zero column in~\eqref{E:13} and rearrange the columns somewhat, we get the C-spline base, $B_{C}$:
\begin{equation}
    \label{E:tja}
B_{C} = \begin{pmatrix} 
 1 &y &x -\Delta x &0\\
 1 &y &0 &x -\Delta x
\end{pmatrix}
\end{equation}
Let
\[
\mathbf{b} =  \begin{pmatrix} b_{1}&b_{2}&b_{3}&b_{4}\end{pmatrix}^{T}
\]
be an arbitrary coefficient vector. Then
\[
\mathbf{z} = B_{C}\mathbf{b}
\]
corresponds with the polynomial equations
\begin{align}
\label{E:17}
z_{1}\left(x,y\right) = b_{1} + b_{2} y + b_{3} \left(x - \Delta x\right), \quad \quad &\;\;\;0\leq x\leq \Delta x,\quad \;\;0\leq y\leq \Delta y \notag \\
z_{2}\left(x,y\right) = b_{1} + b_{2} y + b_{4} \left(x - \Delta x\right), \quad \quad &\Delta x\leq x\leq 2\Delta x,\quad 0\leq y\leq \Delta y
\end{align}
Now, if we substitute $x = \Delta x$ in~\eqref{E:17} we have that for any choice of $\mathbf{b}$ constraint~\eqref{E:8} is satisfied:
\begin{equation}
\label{E:tja2}
z_{1}\left(\Delta x,y\right) = z_{2}\left(\Delta x,y\right) = b_{1} + b_{2} y 
\end{equation}
It follows that $B_{C}$,~\eqref{E:tja}, is the base that enforces zeroth order continuity. 

We summarize, $C^{0}$ continuity between two piecewise polynomials results in a smoothness matrix  $H$,~\eqref{E:9}. The coefficients $\mathbf{c}$,~\eqref{E:6}, defined on the original polynomial base $B$,~\eqref{E:5}, are constrained to lie within the null-space of this smoothness matrix. Stated differently, the coefficients $\mathbf{c}$ are constrained to be a linear combination of the columns of $H_{0}$,~\eqref{E:10}, which span the null space of $H$. By directly multiplying the null-space matrix $H_{0}$ with the the original polynomial base $B$ we get the null-base  $B_{0}$,~\eqref{E:13}, which contains redundant columns consisting of zero vectors. Dropping these zero vectors we obtain the C-spline base $B_{C}$,~\eqref{E:tja}, which has the $C^{0}$ constraint~\eqref{E:8} build into its structure, as may be checked,~\eqref{E:tja2}. 

\section{Enforcing First Order Continuity \label{S:4}}
In order for the partial derivatives of the two polynomials~\eqref{E:3} to connect at the boundary, that is, in order to have $C^{1}$ continuity, we must have that the partial derivatives $\partial z_{1}/\partial x$ and $\partial z_{2}/\partial x$ are $C^{0}$ at their boundaries, that is
\begin{equation}
\label{E:19}
\left.\frac{\partial z_{1}\left(x,y\right)}{\partial x}\right|_{x = \Delta x} = \left.\frac{\partial z_{2}\left(x,y\right)}{\partial x}\right|_{x = \Delta x}
\end{equation}
Substituting~\eqref{E:3} in ~\eqref{E:19}, we find   
\[
c_{12} =  c_{22}
\]
or, equivalently,
\begin{equation}
\label{E:21}
c_{12} -  c_{22} = 0
\end{equation}
Adding constraint~\eqref{E:21} to~\eqref{E:9}, the new smoothness matrix $H$ and the corresponding null space $H_{0}$  become, respectively,
 \[
H =  \begin{pmatrix} 
&1 &\Delta x &y &-1 &0 &-1\\
&0 &1 &0 &0 &-1 &0 \end{pmatrix}
\]
and
\begin{equation}
\label{E:22}
H_{0} =  \begin{pmatrix} 
y &-\Delta x &1 &-y \\
0 &1 &0 &0 \\
0 &0 &0 &1 \\
0 &0 &1 &0 \\
0 &1 &0 &0 \\
1 &0 &0 &0  
\end{pmatrix}
\end{equation}
Multiplying~\eqref{E:22} with the original polynomial base~\eqref{E:5} we get
\begin{align}
B_{0} &= B H_{0}\notag\\
&= \begin{pmatrix} 
y &x -\Delta x 	&1	&0 \\
y &x -\Delta x	&1 	&0
\end{pmatrix}\nonumber
\end{align}
Dropping the redundant zero column and rearranging the columns somewhat, we get the C-spline base:
\begin{equation}
    \label{E:24}
B_{C} = \begin{pmatrix} 
 1 &y &x -\Delta x \\
 1 &y &x -\Delta x
\end{pmatrix}
\end{equation}
which, since $\Delta x$  is a constant, is equivalent to the base
\begin{equation}
    \label{E:25}
B_{C} = \begin{pmatrix} 
 1 &y &x \\
 1 &y &x
\end{pmatrix}
\end{equation}
From C-spline-base~\eqref{E:25} it can be seen that the first order piecewise polynomials which have first order partial derivatives everywhere collapse to a global polynomial of order $d = 1$ and $C^{1}$, which is just the base of a linear regression model having an intercept and predictors $x$ and $y$. 

The here given framework for deriving C-spline bases may be generalized to $d$th order piecewise polynomials with $r$th order continuity, $0\leq r\leq d$, on arbitrary geometries. If one does this then it is found that the C-spline base, $B_{C}$, has a relatively simple structure. This simple structure allows us to directly construct $B_{C}$  without first having to compute the null matrix $H_{0}$ and then taking its matrix product with the original base $B$. This makes C-spline modeling, as given in the next section, computationally efficient. It will be seen that the computational burden of constructing a C-spline is equivalent to that of performing an ordinary regression analysis.

\section{An Algorithm to Construct C-Splines \label{S:5}}
We give here the algorithm for the construction of C-splines for bivariate geometries, partitioned into $I\times J$  adjacent Cartesian domains on which $d$th order piecewise polynomials with $r$th order continuity everywhere are defined.

\subsection{The Geometry}\label{Sub:5}
\noindent First we define a partitioning of the Cartesian $\left(x, y\right)$ plane. Then we number the resulting partitionings. In the region of interest the $x$ values take on values from $a_{x}$ to $b_{x}$ and the $y$ values take on values from $a_{y}$ to $b_{y}$. If we partition the $x$-axis in $I$  adjacent axes with equal lengths $\Delta x$ and the  $y$-axis in $J$  adjacent axes with equal lengths  $\Delta y$. Then this results in $K = IJ$ partitionings. 

Now, we may number each partitioning in the following manner. For $i = 1$ we number the partitionings of the $y$-axis from $k = 1, \ldots, J$ , for $i = 2$ we number the partitionings of the $y$-axis from $k = J + 1, \ldots, 2 J$  , etc… We then have that the $\left(i,j\right)$th partitioning is numbered as
\begin{equation}
\label{E:A2}
k = \left(i - 1\right) J + j,\qquad \qquad \qquad 1\leq i\leq I,\qquad 1\leq j\leq J
\end{equation}
where $k = 1, \ldots, K$ and $K=IJ$. 

In the next paragraph we will construct our C-spline base. The geometry, as given in~\eqref{E:A2}, is non-trivial in that the Cartesian coordinate system having coordinates $\left(i,j\right)$ corresponds with the $k$th row of this C-spline base. 

\subsection{Constructing the C-Spline Base}
\noindent First we construct the building blocks of our base. Let 
 \begin{equation}\label{E:A3}
u_{i} = \begin{cases}
\left(x-a_{x}\right)-i\Delta x, \quad\quad & i = 1,\ldots,I-1\\
\left(x-a_{x}\right)-\left(i-1\right)\Delta x, \quad\quad &i = I
\end{cases}
\end{equation}
Then the $x$-columns of the building blocks are:	
\begin{equation}\label{E:A4}
u_{ki} = \begin{cases}
u_{i}, \quad\quad &k=1,\ldots,iJ, \quad\quad 
\begin{cases}
i = 1,\ldots,I
\end{cases}\\
0, \quad\quad &else
\end{cases}
\end{equation}
where $k = 1, \ldots, K$ and $K=IJ$. Likewise, let 
\begin{equation}\label{E:A5}
v_{j} = \begin{cases}
\left(y-a_{y}\right)-j\Delta y, \quad\quad & j = 1,\ldots,J-1\\
\left(y-a_{y}\right)-\left(j-1\right)\Delta y, \quad\quad &j = J
\end{cases}
\end{equation}
Then the $y$-columns of the building blocks are:		
\begin{equation}\label{E:A6}
v_{kj} = \begin{cases}
v_{j}, \quad\quad &k=1,\ldots,j+\left(i-1\right)\times J, \quad\quad 
\begin{cases}
j = 1,\ldots,J\\
i = 1,\ldots,I
\end{cases}\\
0, \quad\quad &else
\end{cases}
\end{equation}
where $k = 1, \ldots, K$ and $K=IJ$.

Using the building blocks~\eqref{E:A4} and~\eqref{E:A6}, we may now construct the C-spline base $B_{C}$. Our polynomial is of order $d$, that is, let $p$ and $q$ be the powers of $x$ and $y$, respectively, then $0\leq p + q\leq d$. Let
\begin{equation}\label{E:A7}
U_{p} = \begin{cases}
u_{kI}^{p}, \quad\quad & p \leq r\\
\left\{u_{k1}^{p},\ldots,u_{kI}^{p}\right\}, \quad\quad & p > r
\end{cases}
\end{equation}
\\
\begin{equation}\label{E:A8}
V_{q} = \begin{cases}
v_{kJ}^{q}, \quad\quad & q \leq r\\
\left\{v_{k1}^{q},\ldots,v_{kJ}^{q}\right\}, \quad\quad & q > r
\end{cases}
\end{equation}
Then we take the outer product of $U_{p}$ and $V_{q}$ to get $B_{p,q}$, the C-spline equivalent of the polynomial term $x^p y^q$:
\begin{equation}\label{E:A9}
B_{p,q} = U_{p} \otimes V_{q} = \begin{cases}
u_{kI}^{p} v_{kJ}^{q}, \quad\quad & p \leq r, q \leq r\\
\\
\left\{u_{kI}^{p} v_{k1}^{q}, u_{kI}^{p} v_{k2}^{q},\ldots,u_{kI}^{p} v_{kJ}^{q}\right\}, \quad\quad & p \leq r, q > r\\
\\
\left\{u_{k1}^{p} v_{kJ}^{q}, u_{k2}^{p} v_{kJ}^{q},\ldots,u_{kI}^{p} v_{kJ}^{q}\right\}, \quad\quad & p > r, q \leq r\\
\\
\left\{u_{k1}^{p} v_{k1}^{q},u_{k1}^{p} v_{k2}^{q},\ldots,u_{kI}^{p} v_{kJ}^{q}\right\}, \quad\quad & p > r, q > r\\
\end{cases}
\end{equation}
Just as the collection of terms $\left\{x^p y^q\right\}_{0\leq p + q\leq d}$ span the polynomial $P_{d}$,~\eqref{E:0}, So the collection of column vectors
\begin{equation}\label{E:A11}
B_{C}\left(x, y\right) = \left\{B_{p,q}\right\}_{0\leq p + q\leq d}
\end{equation}
span the piecewise polynomials that make up the C-spline. 

Note that for the geometry $I=2$, $J=1$, polynomial order $d=1$ and continuity order $r=0$, the C-spline base~\eqref{E:A11} will differ from~\eqref{E:tja} by one column permutation. Both bases may be considered equivalent in that they both enforce constraint~\eqref{E:8}.

\subsection{Assigning Data Points to the C-Spline Base}
\noindent We have $N$ observed data points in the Cartesian $\left(x, y\right)$-plane that are related to some observed point on the  $z$-axis through the unknown function $f$, that is
\begin{equation}
\label{E:28}
f\left(x_{n}, y_{n}\right) = z_{n},\qquad\qquad  n = 1, \ldots, N.
\end{equation}
By using base~\eqref{E:A11}, we approximate the unknown function $f$ with a collection of piecewise polynomials of degree $d$ that are $C^{r}$ continuous everywhere. To do this we first have to assign each data point $\left(x_{n},y_{n}\right)$ to its corresponding partitioning. 

The $x$- and $y$-axes of each partitioning have,see paragraph~\ref{Sub:5}, lengths of
\[
	\Delta x = \frac{b_{x}-a_{x}}{I},\quad \quad \Delta y = \frac{b_{y}-a_{y}}{J}
\]
We then have that for the data point $\left(x_{n},y_{n}\right)$ which lies in the partitioning having coordinates $\left(i,j\right)$: 
\begin{align}
	a_{x} + \left(i-1\right)\Delta x\leq x_{n}&\leq a_{x} +i\Delta x,  \notag \\ 
	a_{y} + \left(j-1\right)\Delta y\leq y_{n}&\leq a_{y} +j\Delta y  \notag 
\end{align}
or, equivalently,
\[
	\left(i-1\right)\leq \frac{x_{n} - a_{x}}{\Delta x}\leq i,  \quad \quad  \left(j-1\right)\leq \frac{y_{n} - a_{y}}{\Delta y}\leq j
\]
It follows that the coordinates of the partitioning in which the data point $\left(x_{n},y_{n}\right)$ lies may be found as
\begin{equation}
\label{E:29}
	i = ceil\left(\frac{x_{n} - a_{x}}{\Delta x}\right), \quad j = ceil\left(\frac{y_{n} - a_{y}}{\Delta y}\right)
\end{equation}
where $ceil\left( x\right)$ is the function that gives the smallest integer that is greater than or equal to $x$. Substituting these coordinates in~\eqref{E:A2}, we may assign the data point $\left(x_{n},y_{n}\right)$ to its corresponding piecewise polynomial, or, equivalently, to its corresponding row $k$ in the base~\eqref{E:A11}. 
\\
\\
\textit{Example}: \\
Say, we use the C-spline base as given in~\eqref{E:tja}  	
\[
B_{C}\left(x, y\right) = \begin{pmatrix} 
 1 &y &x -\Delta x &0\\
 1 &y &0 &x -\Delta x
\end{pmatrix}
\]
where the first and second row of $B_{C}$  correspond, respectively, with the first and second partitioning of Figure 1. Now, say we have a small dataset of $N = 5$ observations $\left(x_{n},y_{n}\right)$ having values of
\begin{align}
	\left(x_{1},y_{1}\right) &= \left(1.1 \Delta x, 0.3 \Delta y\right) \notag\\
	\left(x_{2},y_{2}\right) &= \left(1.2 \Delta x, 0.7 \Delta y\right) \notag\\
	\left(x_{3},y_{3}\right) &= \left(0.1 \Delta x, 0.3 \Delta y\right) \notag\\
	\left(x_{4},y_{4}\right) &= \left(0.5 \Delta x, 0.1 \Delta y\right) \notag\\
	\left(x_{5},y_{5}\right) &= \left(1.7 \Delta x, 0.8 \Delta y\right) \notag
\end{align}
where $\Delta x$ and $\Delta y$ are some constants. Then, using~\eqref{E:A2} and~\eqref{E:29}, the points $\left(x_{3},y_{3}\right)$ and  $\left(x_{4},y_{4}\right)$  are assigned to the first partitioning, or, equivalently, to the first row of $B_{C}$. Likewise, $\left(x_{1},y_{1}\right)$, $\left(x_{2},y_{2}\right)$ and  $\left(x_{5},y_{5}\right)$ are assigned to the second partitioning, or, equivalently, to the second row of $B_{C}$:
	\begin{equation}
	\label{E:29b}
	\tilde{B}_{C} = \begin{pmatrix} 
 B_{C}^{\left(2\right)}\left(x_{1}, y_{1}\right)\\
 B_{C}^{\left(2\right)}\left(x_{2}, y_{2}\right)\\
 B_{C}^{\left(1\right)}\left(x_{3}, y_{3}\right)\\ 
 B_{C}^{\left(1\right)}\left(x_{4}, y_{4}\right)\\
 B_{C}^{\left(2\right)}\left(x_{5}, y_{5}\right)
\end{pmatrix} = \begin{pmatrix} 
 1 &0.3\Delta y &0 							&0.1\Delta x \\
 1 &0.7\Delta y &0 							&0.2\Delta x \\
 1 &0.3\Delta y &-0.9\Delta x 	&0 		\\ 
 1 &0.1\Delta y &-0.5\Delta x 	&0 		\\
 1 &0.8\Delta y &0 							&0.7\Delta x
\end{pmatrix}
	\end{equation}
Note that we use a tilde to signify a base $B_{C}$ to which data points have been assigned.

\subsection{Constructing a C-Spline}
Let $m$ be the number of columns of the C-spline base $B_{C}$,~\eqref{E:A11}. Then, after we have assigned all $N$ data points to the base $B_{C}$, we get the $N\times m$ matrix $\tilde{B}_{C}$, see~\eqref{E:29b}. The unknown $\textbf{b}$ coefficients, see~\eqref{E:17}, of the C-spline may be found as the least-squares solution of 
\begin{equation}
\label{E:30}
	\textbf{b} = \left(\tilde{B}_{C}^{T}\tilde{B}_{C}\right)^{-1}\tilde{B}_{C}^{T}\textbf{z}
\end{equation}
where $\textbf{z} =  \begin{pmatrix} z_{1} &\cdots &z_{N}\end{pmatrix}$ is the vector with output values,~\eqref{E:28}.

Now, say we wish to get the C-spline estimate $\hat{z}_{N+1}$ of the data point $\left(x_{N+1}, y_{N+1}\right)$. Then, using~\eqref{E:A2} and~\eqref{E:29}, we first determine the row $k$ of the base $B_{C}$,~\eqref{E:A11}, that corresponds with this data point and then plug in its value. This results in the $1\times m$ row-vector
\begin{equation}
\label{E:30a}
\tilde{B}_{C}^{\left(k\right)} = B_{C}^{\left(k\right)}\left(x_{N+1},y_{N+1}\right)
\end{equation}
The estimate $\hat{z}_{N+1}$ is then found by simply taking the inner product of~\eqref{E:30a} and~\eqref{E:30}:
\[
\hat{\textbf{z}} = \tilde{B}_{C}^{\left(k\right)}\cdot \textbf{b}
\]
We see that constructing a C-spline is equivalent to performing a regression analysis.

\section{Discussion \label{S:9}}
We have introduced here Cartesian splines, or C-splines, for short. C-splines are piecewise polynomials which are defined on  adjacent Cartesian coordinate systems and are $C^{r}$ continuous throughout. We have given here an algorithm that allows one to construct C-spline bases without first having to find the null-space of the corresponding smoothness matrix $H$. 

This makes the construction of a given C-spline base computationally trivial since no null-space of $H$ has to be evaluated. This means that for C-splines the computational burden lies solely, just as in any ordinary regression analysis, in the evaluation of the inverse of $\tilde{B}_{C}^{T}\tilde{B}_{C}$, where $\tilde{B}_{C}$ is the matrix with the independent variables. 

Note that the algorithm, equations~\eqref{E:A2} through~\eqref{E:A11}, may be generalized relatively easy to construct C-splines for multivariate domains.
\\
\\
\noindent\textbf{Acknowledgments:} This research was partly funded by the Delft Cluster project (no. CT04.33.11), which is kindly acknowledged by the authors.

\end{document}